\documentclass[reqno,11pt]{amsart}
\usepackage[margin=4cm]{geometry}
\usepackage{amscd,amssymb}
\usepackage{mathtools}
\usepackage{color}
\usepackage{hyperref}
\usepackage{comment}
\usepackage{tikz}
\usepackage{slashed}
\usetikzlibrary{matrix}
\usepackage{marginnote}
\usepackage{enumitem}
\usepackage{amsmath,amsfonts,latexsym}

\numberwithin{equation}{section}

\theoremstyle{plain}
\newtheorem{theorem}[subsection]{Theorem}
\newtheorem{lemma}[subsection]{Lemma}
\newtheorem{corollary}[subsection]{Corollary}

\theoremstyle{definition}

\theoremstyle{remark}
\newtheorem{remark}[subsection]{Remark}

\newcommand{\RR}{\mathbb{R}}

\usepackage{anysize}

\begin{document}

\title{A note on the Gaussian curvature on noncompact surfaces.}

\author{Simone Cecchini}
\address{Department of Mathematics,
Northeastern University,
Boston, MA 02115,
USA}

\email{cecchini.s@husky.neu.edu}


\begin{abstract} 
We give a short proof of the following fact.
Let $\Sigma$ be a connected, finitely connected noncompact manifold without boundary. 
If $g$ is a complete Riemannian metric on $\Sigma$ whose Gaussian curvature $K$ is nonnegative at infinity, then $K$ is integrable.
In particular, we obtain a new short proof of the fact if $\Sigma$ admits a complete metric whose Gaussian curvature is nonnegative and positive at one point, then $\Sigma$ is diffeomorphic to $\RR^2$.
\end{abstract}


\maketitle

\section{Introduction}
In~\cite{RosSto94}, J.Rosenberg and S.Stolz conjectured that a closed manifold $X$ admits a metric of positive scalar curvature when the cylinder $X\times\RR$ admits a \emph{complete} metric of positive scalar curvature.
When $X$ is one-dimensional, this conjecture corresponds to the statement that the cylinder $S^1\times\RR$ cannot carry a complete Riemannian metric of positive scalar curvature.
This fact is well-known and follows for instance from~\cite[Corollary~6.13]{GroLaw83}.
One of the main aims of this note is to provide an elementary proof of this fact.

More generally, we consider a class of noncompact surfaces, which are called finitely connected.
Roughly speaking, they are obtained by removing a finite number of points from a compact surface.
For such surfaces, the Euler characteristic is obviously well-defined.
Suppose that $g$ is a complete Riemannian metric on a connected, finitely connected, noncompact surface without boundary $\Sigma$.
We show that if the Gaussian curvature $K$ of $g$ is nonnegative at infinity, then $K$ must be integrable and its integral must be smaller than the Euler characteristic of $\Sigma$.
As an application, we obtain a new proof of the fact that a connected, finitely connected noncompact surface without boundary with nonpositive Euler characteristic cannot carry a complete metric whose Gaussian curvature is nonnegative and positive at one point.
In particular, this applies to the cylinder $S^1\times\RR$, since it is finitely connected and $\chi(S^1\times\RR)=\chi(S^1)=0$. 


\section{The theorem}


Before stating the theorem, we review some basic notions on noncompact surfaces.
We say that a smooth surface without boundary $\Sigma$ is \emph{finitely connected} if there exists a compact surface with boundary $\Omega\subset \Sigma$ such that
\begin{enumerate}[label=(\roman*)]
\item the boundary of $\Omega$ is a disjoint union of closed simple curves $l_1,\ldots, l_p$;  
\item the open set $\Sigma\setminus \Omega$ is a disjoint union of the cylinders $C_j:=l_j\times (0,\infty)$, for $j=1,\ldots p$.
\end{enumerate}
Notice that equivalently a noncompact surface without boundary $\Sigma$ is finitely connected if it is homeomorphic to a closed surface $\Sigma_1$ with $p$ points removed.
In this case, the singular homology groups of $\Sigma$ have finite rank and the Euler characteristic $\chi(\Sigma)$ of $\Sigma$ is given by the formula
\begin{equation}\label{E:Euler char. finitely connected}
	\chi(\Sigma)\ =\ \chi(\Sigma_1)\,-\,p\,.
\end{equation}
For the basic notions on finitely connected surfaces, we refer to~\cite[Section~2.1]{ShiShiTan03}.


Let $(\Sigma,g)$ be a noncompact Riemannian surface .
Let $K$ be the Gaussian curvature of $g$, and $dA_g$ the area element.
We say that a Riemannian surface $(\Sigma,g)$ admits \emph{total curvature} if 
\[
	\int_\Sigma K_+\,dA_g\,<\infty\quad\qquad\text{ or }\quad\qquad\int_\Sigma K_-\,dA_g\,<\infty\,,
\]
where $K_+=\max\{K,0\}$ and $K_-:=\max(-K,0)$.
In this case, the extended real number
\begin{equation}\label{E:total curvature}
	c(\Sigma;g)\ :=\ \int_\Sigma K\,dA_g\ =\ 
	\int_\Sigma K_+\,dA_g\,-\,\int_\Sigma K_-\,dA_g\in[-\infty,\infty]
\end{equation}
is called the \emph{total curvature} of $(\Sigma,g)$ (for more details on the notion of total curvature for noncompact surfaces, see~\cite[Section~2.1]{ShiShiTan03}).
We say that the Gaussian curvature $K$ is \emph{nonnegative at infinity} if $K\geq 0$ outside a compact subset of $\Sigma$.
Notice that in this case $K_-$ is compactly supported so that $(\Sigma,g)$ admits total curvature $c(\Sigma;g)$ ranging over the interval $(-\infty,+\infty]$.


\begin{theorem}\label{T:main theorem}
Let $(\Sigma,g)$ be a connected, finitely connected, complete Riemannian surface without boundary.
If $K$ is nonnegative at infinity, then $(\Sigma,g)$ admits \emph{finite} total curtvature.
Moreover, we have
\begin{equation}\label{E:G-B inequality}
	2\pi\,\chi(\Sigma)\ \geq\ c(\Sigma;g)\,.
\end{equation}
\end{theorem}


\begin{remark}
Formula~\eqref{E:G-B inequality} can be deduced from the first part of the theorem by using the Gauss-Bonnet inequality due to Cohn-Vossen (cf.~\cite{Coh35} and~\cite[Theorem~2.2.1]{ShiShiTan03}).
In our proof of Theorem~\ref{T:main theorem}, we also prove inequality~\eqref{E:G-B inequality} in the case when the Gaussian curvature is nonnegative at infinity.
\end{remark}


Notice that, by formula~\eqref{E:Euler char. finitely connected}, if a connected, finitely connected, noncompact surface has positive Euler characteristic, then it is homeomorphic (and hence diffeomorphic) to $\RR^2$.
Therefore, from Theorem~\ref{T:main theorem} we obtain a new proof of the following fact.


\begin{corollary}\label{C:Cheeger-Gromoll}
Let $\Sigma$ be a connected, finitely connected, noncompact surface.
If $\Sigma$ admits a complete metric whose Gaussian curvature is nonnegative and positive at one point, then $\Sigma$ is diffeomorphic to $\RR^2$.
\end{corollary}


\begin{remark}
Corollary~\ref{C:Cheeger-Gromoll} is a direct consequence of the soul conjecture of Cheeger and Gromoll, proved in full generality by Perelman in~\cite{Per94}, when specialized to the two-dimensional case.
\end{remark}


\section{The approximation procedure}
In this section we present the proof of Theorem~\ref{T:main theorem}.

Let $\Sigma$ be a noncompact, connected, finitely connected surface without boundary.
Let $\Omega\subset \Sigma$ be a compact submanifold with boundary such that the boundary $\partial\Omega$ of $\Omega$ consists of $p$-copies of $S^1$ and $\Sigma\setminus\Omega=\bigsqcup_{j=1}^pC_j$, where each $C_j$ is a copy of the cylinder $S^1\times (0,\infty)$.
We also assume $\Sigma$ is \emph{orientable} and pick an orientation.

Let $\Sigma_h$ be the compact surface with boundary obtained by truncating the cylindrical ends of $\Sigma$ at the height $h$.
This means that the boundary $\partial\Sigma_h$ of $\Sigma_h$ is the disjoint union of $p$ copies of $S^1$ and $\Sigma\setminus\Sigma_h=\bigsqcup_{j=1}^p\{S^1\times (h,\infty)\}$.
The \emph{total geodesic curvature} of $\Sigma_h$ is defined by
\begin{equation}\label{E:lambda}
	\lambda(h)\ :=\ \int_{\partial\Sigma_h}K\,,
\end{equation}
where the boundary $\partial\Sigma_h$ is \emph{positively} oriented with respect to the given orientation of $\Sigma$ and $K$ is the \emph{geodesic curvature} of $\partial\Sigma_h$ (see~\cite[Section~4-4, Definition~10]{docarmo}).


\begin{lemma}\label{L:lambda>=0}
Let $(\Sigma,g)$ be a connected, finitely connected, \emph{orientable}, complete Riemannian surface without boundary.
If the Gaussian curvature of $g$ is nonnegative at infinity, then $\lambda(h)$ converges to a nonnegative number $L$, as $h$ goes to infinity.
\end{lemma}


\begin{remark} 
Our proof of this lemma is based on the fact that we can choose on each cylindrical end $C_j$ suitable coordinates which simplify the components of the metric.
Such coordinates were first used by S. Rosenberg in~\cite{SRos82} to provide a short proof of Cohn-Vossen inequality.
\end{remark}


\subsection{Proof of Lemma~\ref{L:lambda>=0}}
Since $\Sigma_h$ is a retract of $\Sigma$, using the Gauss-Bonnet theorem (see~\cite[Section~4-5]{docarmo}) on $\Sigma_h$ we obtain
\begin{equation}\label{E:Gauss-Bonnet}
	2\pi\,\chi (\Sigma)\ =\ 2\pi\,\chi (\Sigma_h)\ =\ c(\Sigma_h;g)\,+\,\lambda (h)\,.
\end{equation}
Since $K$ is nonnegative at infinity, there exists $h_1>0$ such that $K\geq 0$ on $\Sigma\setminus\Sigma_{h}$, for all $h\geq h_1$.
From~\eqref{E:Gauss-Bonnet} we deduce that the function
\begin{equation}\label{E:L1}
	\lambda(h)\ =\ 2\pi\,\chi(\Sigma)\,-\,c(\Sigma_h;g)
\end{equation}
is nonincreasing on the interval $(h_1,\infty)$, so that the extended real number 
\begin{equation}\label{E:L3}
	L\ :=\ \lim_{h\rightarrow\infty}\lambda(h)\in [-\infty,+\infty)
\end{equation}
is well-defined.
To conclude the proof, it remains to show that we must have $L\geq 0$.

In order to get more information on the function $\lambda(h)$ and the number $L$, we compute a local expression for the geodesic curvature $K$.
As observed in~\cite{SRos82}, we can choose coordinates $(t_j,\theta_j)$ on the cylindrical end $C_j$ in a way that 
\begin{enumerate}[label=(\roman*)]
	\item for all $P\in C_j$, the basis
	\[
		\left\{\frac{\partial}{\partial t_j}\Big|_P \,,\   \frac{\partial}{\partial\theta_j}\Big|_P\right\}
	\] 
	of the tangent space $T_PC_j\cong \RR\oplus\RR$ is positively oriented;
	\item the metric $g$, restricted to the cylindrical end $C_j$, is of the form
	\[
		g(t_j,\theta_j)\ =\ dt_j^2\,+\,G_j(t_j,\theta_j)\,d\theta_j^2\,,\qquad\qquad (t_j,\theta_j)\in C_j\,,
	\]
	where $G_j:C_j\rightarrow (0,\infty)$ is a smooth function.
\end{enumerate}
For the rigorous construction of the coordinates $(t_j,\theta_j)$, we refer to~\cite[page 747]{SRos85}.
With this choice, the curve $\gamma_j^h(s)=(t_j(s),\theta_j(s))=(h,s)$ parametrizes $\partial \Sigma_h\cap C_j$ with \emph{positive} orientation  (see~\cite[pp.~267-268]{docarmo}).
Moreover, by~\cite[Section~4-4, Proposition~3]{docarmo}, the geodesic curvature $K$ of $\gamma_j^h$ takes the form 
\begin{equation}\label{E:geodesic curvature}
	K\ =\ \frac{1}{2\,\sqrt{G_j}\,}\,\frac{\partial G_j}{\partial t_j}\,\frac{d\theta_j}{ds}
	\ =\ \frac{\partial}{\partial t_j}\,\sqrt{G_j}\,,
\end{equation}
from which
\[
	\int_{\partial\Sigma_h\cap C_j}K\ =\ 
	\int_0^{2\pi}\left(\frac{\partial}{\partial t_j}\,\sqrt{G_j}\,\right)(h,s)\,ds\ =\ 
	\frac{d}{dh}\,\int_0^{2\pi}\sqrt{G_j(h,s)}\,ds\,.
\]
Hence,
\begin{equation}\label{E:L2}
	\lambda(h)\ =\ \mu^\prime(h)\,,
\end{equation}
where
\begin{equation}\label{E:mu}
	\mu(h)\ :=\ \sum_{j=1}^p\,\int_0^{2\pi}\sqrt{G_j(s,h)}\,ds\,.
\end{equation}

Finally, we use Equation~\eqref{E:L2} to deduce that the number $L$ defined by~\eqref{E:L3} must be nonnegative.
Suppose indeed that $L<0$. 
Then by Equation~\eqref{E:L2} it follows that $\lim_{h\rightarrow \infty}\mu(h)=-\infty$, which is impossible, since, by the expression~\eqref{E:mu}, $\mu(h)$ is a strictly positive function.
Therefore, we must have $L\geq 0$, which concludes the proof.
\hfill$\square$


\subsection{Proof of Theorem~\ref{T:main theorem}}
Let $\Sigma$, $g$ be as in the hypothesis of Theorem~\ref{T:main theorem}.
We also assume $\Sigma$ is \emph{orientable}: the case when $\Sigma$ is nonorientable is obtained, by a standard argument, by considering the orientable double cover.
From~\eqref{E:Gauss-Bonnet}, we have
\[
	c(\Sigma_h;g)\ =\ 2\pi\,\chi(\Sigma)\,-\,\lambda(h)\,.
\]
Taking the limit for $h\rightarrow\infty$ and using Lemma~\ref{L:lambda>=0}, we deduce that
\begin{equation}\label{E:G-B nonnegative at infty}
	c(\Sigma;g)\ =\ 2\pi\,\chi(\Sigma)\,-\,L\ \leq\ 2\pi\,\chi(\Sigma)\,.
\end{equation}
Therefore, $c(\Sigma;g)$ is finite and satisfies the Gauss-Bonnet Inequality~\eqref{E:G-B inequality}.
\hfill$\square$


\end{document}